\documentclass[12pt]{amsart}
\usepackage[T1]{fontenc}
\usepackage[latin1]{inputenc}
\usepackage{times}
\makeatletter

\newcommand{\LyX}{L\kern-.1667em\lower.25em\hbox{Y}\kern-.125emX\spacefactor1000}

\theoremstyle{plain}    
\newtheorem{thm}{Theorem}[section]
\numberwithin{equation}{section} %% Comment out for sequentially-numbered
\numberwithin{figure}{section} %% Comment out for sequentially-numbered
\theoremstyle{plain}    
\newtheorem*{thm*}{Theorem} 
\theoremstyle{plain}    
\newtheorem{cor}[thm]{Corollary} %%Delete [thm] to re-start numbering
\theoremstyle{plain}    
\newtheorem{prop}[thm]{Proposition} %%Delete [thm] to re-start numbering

\makeatother

\begin{document}

\title[Characterization of freeness]{A Characterization of Freeness \\
by a \\
Factorization Property of \protect\( R\protect \)-transform.}

\author[A. Nica]{Alexandru Nica{*} }

\address{Department of Pure Mathematics, University of Waterloo, Waterloo, Ontario,
N2L 3G1, Canada}

\email{anica@math.uwaterloo.ca}

\thanks{{*} Research supported by a grant of NSERC, Canada.}

\author[D. Shlyakhtenko]{Dimitri Shlyakhtenko\dag}

\address{Department of Mathematics, UCLA, Los Angeles, CA 90095, USA}

\email{shlyakht@math.ucla.edu}

\thanks{\dag Research supported by an NSF postdoctoral fellowship.}

\author[R. Speicher]{Roland Speicher}

\address{Department of Mathematics and Statistics, Queen's University, Kingston, Ontario
K7L 3N6, Canada}

\email{speicher@mast.queensu.ca}

\begin{abstract}
Let \( M \) be a \( B \)-probability space. Assume that \( B \) itself is
a \( D \)-probability space; then \( M \) can be viewed as a \( D \)-probability
space as well. Let \( X\in M \). We characterize freeness of \( X \) from
\( B \) with amalgamation over \( D \) in terms of a certain factorization
condition linking the \( B \)-valued and \( D \)-valued \( R \)-transforms
of \( X \). We give an application to random matrices.
\end{abstract}
\maketitle

\section{Operator-valued \protect\( R\protect \)-transform.}

\subsection{Multiplicative functions.}

Let \( B \) be a unital algebra. Recall that a \emph{\( B \)-probability space}
\( (M,E:M\to B) \) (see e.g. \cite{DVV:book}, \cite{speicher:thesis}) is
a pair consisting of an algebra \( M \) containing \( B \) as a unital subalgebra,
and a \emph{conditional expectation}~\( E:M\to B \). In other words, \( E:M\to B \)
is unital and \( B \)-bilinear:
\[
E(1)=1,\qquad E(bmb')=bE(m)b',\quad \forall b,b'\in B,\forall m\in M.\]
Elements of \( M \) are called \( B \)-valued random variables. 

Recall that a \emph{multiplicative map} \( \langle \cdots \rangle :\bigcup M^{n}\to B \)
is a \( \mathbb {C} \)-multilinear map, (i.e., a sequence of maps \( M^{n}\ni m_{1},\dots ,m_{n}\mapsto \langle m_{1},\dots ,m_{n}\rangle \in B \))
satisfying the \( B \)-linearity conditions
\begin{eqnarray*}
\langle bm_{1},\dots ,m_{n}\rangle  & = & b\langle m_{1},\dots ,m_{n}\rangle \\
\langle m_{1},\dots ,m_{n}b\rangle  & = & \langle m_{1},\dots ,m_{n}\rangle b\\
\langle m_{1},\dots ,m_{k}b,m_{k+1},\dots ,m_{n}\rangle  & = & \langle m_{1},\dots ,m_{k},bm_{k+1},\dots ,m_{n}\rangle .
\end{eqnarray*}
(Here \( m_{i}\in M \), \( b\in B \)). 

Given a non-crossing partition \( \pi \in NC(n) \) and an arbitrary multiplicative
map \( \langle \cdots \rangle  \), we can construct a \emph{bracketing}~\( \pi \langle \cdots \rangle :M^{n}\to B \),
defined recursively by
\begin{eqnarray*}
1_{k}\langle m_{1},\dots ,m_{k}\rangle  & = & \langle m_{1},\dots ,m_{k}\rangle \\
(\pi \sqcup \rho )\langle m_{1},\dots ,m_{k}\rangle  & = & \pi \langle m_{1},\dots ,m_{p}\rangle \cdot \rho \langle m_{p+1},\dots ,m_{k}\rangle \\
\textrm{ins}(p,\rho \to \pi )(m_{1},\dots ,m_{k}) & = & \pi \langle m_{1},\dots ,m_{p}\rho \langle m_{p+1},\dots ,m_{p+q}\rangle ,\\
 &  & \quad m_{p+q+1},\dots ,m_{k}\rangle .
\end{eqnarray*}
Here \( 1_{k} \) denotes the partition with the sole class \( \{1,\dots ,k\} \),
\( \pi \sqcup \rho  \) denotes disjoint union (with the equivalence classes
of \( \rho  \) placed after those of \( \pi  \)), and \( \textrm{ins}(p,\rho \to \pi ) \)
denotes the partition obtained from \( \pi  \) by inserting the partition \( \rho  \)
after the \( p \)-th element of the set on which \( \pi  \) determines a partition. 

In other words, each partition \( \pi  \) is interpreted as a recipe for placing
brackets \( \langle \cdots \rangle  \), and \( \pi \langle \cdots \rangle  \)
is the value of the resulting expression.

\subsection{Moments and \protect\( R\protect \)-transform.}

The \( B \)-probability space structure of the algebra \( M \) gives rise
to one example of such a multiplicative map, namely, the \emph{moments map}
\[
[\cdots ]_{E}:\bigcup M^{n}\to B,\]
given by
\[
[m_{1},\dots ,m_{n}]_{E}=E(m_{1}\cdots m_{n}).\]
The reason for the name is that, having fixed \( B \)-random variables \( X_{1},\dots ,X_{p}\in M \),
the following values of \( [\cdots ]_{E} \),
\[
[b_{0}X_{i_{1}}b_{1},\cdots ,X_{i_{n-1}}b_{n-1},X_{i_{n}}b_{n}]_{E}\]
are called \emph{\( B \)-valued moments} of the family \( X_{1},\dots ,X_{n} \). 

In \cite{speicher:thesis} and \cite{dvv:amalg} the notion of \( B \)-valued
\( R \)-transform was introduced (we follow the combinatorial approach of \cite{speicher:thesis}, see also \cite{cumulants}).
Like the map \( [\cdots ] \), the \( R \)-transform map is a multiplicative
map
\[
\{\cdots \}_{E}:\bigcup M^{n}\to B.\]
 The following combinatorial formula actually determines \( \{\cdots \}_{E} \)
uniquely:
\begin{eqnarray*}
[m_{1},\dots ,m_{n}]_{E}=\sum \textrm{all possible bracketings involving }\{\cdots \}_{E} &  & \\
=\sum _{\pi \in NC(n)}\pi \{m_{1},\dots ,m_{n}\}_{E}. &  & 
\end{eqnarray*}
The uniqueness of the definition can be easily seen by observing that the right-hand
side of the equation above involves \( \{m_{1},\dots ,m_{n}\}_{E} \) and that
the rest of the terms are products of factors of smaller order (i.e., restrictions
of \( \{\cdots \}_{E} \) to \( M^{k} \), \( k<n \)). 

It is important to note that \( [\cdots ] \) determines \( \{\cdots \} \)
and vice-versa. Moreover, the value of \( [\cdots ]|_{M^{n}} \) depends only
on \( \{\cdots \}|_{M\cup \cdots \cup M^{n}} \), and vice-versa.

\subsection{Moment and cumulant series.}

Given a multiplicative function \( \langle \cdots \rangle :\bigcup M^{k}\to B \),
and \( X_{1},\dots ,X_{n}\in B \), consider the family of multilinear maps
\[
M^{\langle \cdots \rangle }_{i_{1},\dots ,i_{k}}:B^{k-1}\to B,\quad i_{1},\dots ,i_{k}\in \{1,\dots ,n\}\]
given by
\[
M^{\langle \cdots \rangle }_{i_{1},\dots ,i_{k}}(b_{1},\dots ,b_{k-1})=\langle X_{i_{1}}b_{1},\cdots ,b_{k-1}X_{i_{k}}\rangle .\]
In the particular examples above, we get the \emph{moment series} of \( X_{1},\dots ,X_{n} \),
\[
\mu ^{X_{1},\dots ,X_{n}}_{i_{1},\dots ,i_{k}}=M^{[\cdots ]_{E}}_{i_{1},\dots ,i_{k}}\]
and the \emph{cumulant series},
\[
k_{i_{1},\dots ,i_{k}}^{X_{1},\dots ,X_{n}}=M_{i_{1},\dots ,i_{k}}^{\{\cdots \}_{E}}.\]
We will sometimes write
\[
k_{B;i_{1},\dots ,i_{k}}^{X_{1},\dots ,X_{n}}\]
to emphasize that the series is valued in \( B \).

\subsection{Freeness with amalgamation.}

Let \( M_{1},M_{2}\subset M \) be two subalgebras, each containing \( B \).
The importance of \( R \)-transform is apparent from the following theorem
(\cite{speicher:thesis}, see also \cite{cumulants-nica}):

\begin{thm*}
Let \( S_{1},S_{2} \) be two subsets of \( M \). Let \( M_{i} \) be the algebra
generated by \( S_{i} \) and \( B \), \( i=1,2. \) Then \( M_{1} \) and
\( M_{2} \) are free with amalgamation over \( B \) iff whenever \( X_{1},\dots ,X_{n}\in S_{1}\cup S_{2} \),
\[
\{X_{1},\dots ,X_{n}\}_{E}=0\]
unless either all \( X_{1},\dots ,X_{n}\in S_{1} \), or all \( X_{1},\dots ,X_{n}\in S_{2} \).
\end{thm*}

\subsection{Canonical random variables.\label{sec:canonical}}

Let \( X_{1},\dots ,X_{n}\in M \) be fixed. 

Then by \cite{dvv:amalg} there exists a \( B \)-probability space \( (\mathcal{F},E_{B}:\mathcal{F}\to B) \),
elements \( \lambda _{1}^{*},\dots ,\lambda _{n}^{*}\in \mathcal{F} \), and
elements \( \lambda _{p}^{k} \) satisfying the following properties:\renewcommand{\labelenumi}{(\roman{enumi})}\renewcommand{\theenumi}{\arabic{section}.\arabic{subsection}(\roman{enumi})}

\begin{enumerate}
\item \label{relation:lambdas}\( \lambda _{j_{1}}^{*}b_{1}\lambda _{j_{2}}^{*}b_{2}\dots \lambda _{j_{k}}^{*}b_{k}\lambda _{j}^{k}=k_{j_{1},\dots ,j_{k},j}^{X_{1},\dots ,X_{n}}(b_{1},\dots ,b_{k}) \),
\( b_{1},\dots ,b_{k}\in B \); 
\item \label{relation:E}Let \( w=b_{0}a_{1}b_{1}a_{2}b_{2}\dots a_{n}b_{n} \), where
\( b_{i}\in B \), \( a_{i}=\lambda ^{k}_{p} \) or \( a_{i}=\lambda _{p}^{*} \).
Then \( E_{B}(w)=0 \) unless \( w \) can be reduced to an element of \( B \)
using relation \ref{relation:lambdas}.
\end{enumerate}
(By convention, we take \( \lambda _{p}^{0}=E_{B}(X_{p})\in B \)). We should
clarify that for each \( k \), \( \lambda _{p}^{k} \) is just a formal variable,
and we do not assume any relations between \( \{\lambda _{k}^{p}\}_{k,p} \):
for example, \( \lambda _{p}^{k} \) is not the \( k \)-th power of \( \lambda _{p}^{1} \).

It is not hard to show that the properties listed above determine the restriction
of \( E_{B} \) to the algebra generated by \( \lambda ^{*} \) and \( \{\lambda _{p}^{k}\}_{p,k} \). 

Let
\[
Y_{j}=\lambda _{j}^{*}+\sum _{k\geq 0}\lambda _{j}^{k}.\]
(This series is formal; however, \( Y_{1},\dots ,Y_{n} \) have moments, since
each such moment involves only a finite number of terms from the series defining
\( Y_{j} \)).

It turns out that
\[
k_{i_{1},\dots ,i_{k}}^{X_{1},\dots ,X_{n}}=k_{i_{1},\dots ,i_{k}}^{Y_{1},\dots ,Y_{n}},\quad \textrm{and }\mu _{i_{1},\dots ,i_{k}}^{X_{1},\dots ,X_{n}}=\mu _{i_{1},\dots ,i_{k}}^{Y_{1},\dots ,Y_{n}}.\]
In other words, given a cumulant series, \( Y_{1},\dots ,Y_{n} \) is an explicit
family of \( B \)-valued random variables, whose cumulant series is equal to
the one given.

\section{Freeness from a subalgebra.}

\subsection{\protect\( D\protect \)-cumulants vs. \protect\( B\protect \)-cumulants.}

Let now \( D\subset B \) be a unital subalgebra, and let \( F:B\to D \) be
a conditional expectation. If \( (M,E:M\to B) \) is a \( B \)-probability
space, then \( (M,F\circ E:M\to D) \) is a \( D \)-probability space. 

\begin{thm}
Let \( X_{1},\dots ,X_{n}\in M \). Assume that the \( B \)-valued cumulants
of \( X_{1},\dots ,X_{n} \) satisfy
\[
k^{X_{1},\dots ,X_{n}}_{B;i_{1},\dots ,i_{k}}(d_{1},\dots ,d_{k-1})\in D,\quad \forall d_{1},\dots ,d_{k-1}\in D.\]
Then the \( D \)-valued cumulants of \( X_{1},\dots ,X_{n} \) are given by
the restrictions of the \( B \)-valued cumulants:
\[
k^{X_{1},\dots ,X_{n}}_{D;i_{1},\dots ,i_{k}}(d_{1},\dots ,d_{k-1})=k^{X_{1},\dots ,X_{n}}_{B;i_{1},\dots ,i_{k}}(d_{1},\dots ,d_{k-1}),\quad \forall d_{1},\dots ,d_{k-1}\in D.\]

\end{thm}
\begin{proof}
Let \( N \) be the algebra generated by \( D \) and \( X_{1},\dots ,X_{n} \).
The condition on cumulants implies that
\[
\{\cdots \}_{E}|_{\bigcup N^{p}}\]
is valued in \( D \). It follows from the moment-cumulant formula that \( [\cdots ]_{E}|_{\bigcup N^{p}} \)
is valued in \( D \), and hence that
\[
[\cdots ]_{F\circ E}|_{\bigcup N^{p}}=[\cdots ]_{E}|_{\bigcup N^{p}}=\sum _{\pi \in NC(p)}\pi \{\cdots \}_{E}|_{\bigcup N^{p}}.\]
Since the moment-cumulant formula determines \( \{\cdots \}_{F\circ E}|_{\bigcup N^{p}} \),
it follows that
\[
\{\cdots \}_{F\circ E}|_{\bigcup N^{p}}=\{\cdots \}_{E}|_{\bigcup N^{p}}.\]

\end{proof}
We record an equivalent formulation of the theorem above (which was implicit
in the proof):

\begin{thm}
\label{thrm:restricRtransformDvalued}Let \( N\subset M \) be a subalgebra,
containing \( D \). Assume that \( \{\cdots \}_{E}|_{\bigcup N^{p}} \) is
valued in \( D \). Then
\[
\{\cdots \}_{F\circ E}|_{\bigcup N^{p}}=\{\cdots \}_{E}|_{\bigcup N^{p}}.\]

\end{thm}
In general, in the absence of the condition that \( k^{X_{1},\dots ,X_{n}}_{B} \)
restricted to \( \cup D^{p} \) is valued in \( D \), the expression of \( D \)-valued
cumulants of \( X_{1},\dots ,X_{n} \) in terms of the \( B \)-valued cumulants
is quite complicated. Note, for example, that if \( X \) is a \( B \)-valued
random variable, and \( b\in B \), then the \( B \)-valued cumulant series
of \( bX \) are very easy to describe. On the other hand, the \( D \)-valued
cumulant series of \( bX \) can have a very complicated expression in terms
of the \( D \)-valued cumulant series of \( X \) and \( b \).

The sufficient condition in the theorem above is actually quite close to being
necessary in the case that the conditional expectations are positive maps of
\( * \)-algebras. As an illustration, consider the case that \( D\subset B \)
consists of scalar multiples of \( 1 \), and \( F:B\to D \) is such that \( \tau =F\circ E \)
is a \emph{trace} on \( M \), satisfying \( \tau (xy)=\tau (yx) \) for all
\( x,y\in M \).

Recall that \( X \) is called a \( B \)-semicircular variable if its cumulant
series is given by
\[
k^{X}_{\underbrace{1,\dots 1}_{p}}(b_{1},\dots ,b_{p-1})=\delta _{p,2}\eta (b_{1})\]
for some map \( \eta :B\to B \). It is easily seen that if \( X \) is \( B \)-semicircular,
then \( \eta (b)=E(XbX) \).

\begin{thm}
\label{thrm:characteriRestrictedSemicircular}Let \( (M,E:M\to B) \) be a \( B \)-probability
space, such that \( M \) and \( B \) are \( C^{*} \)-algebras. Let \( F:B\to \mathbb {C}=D\subset B \)
be a faithful state. Assume that \( \tau =F\circ E \) is a faithful trace on
\( M \). Let \( X \) be a \( B \)-semicircular variable in \( M \). Then
the distribution of \( X \) with respect to \( \tau  \) is the semicircle
law iff \( E(X^{2})\in \mathbb {C} \).
\end{thm}
\begin{proof}
If \( E(X^{2})=k_{11}^{X}(1)\in \mathbb {C} \), it follows that the \( B \)-valued
cumulants of \( X \), restricted to \( D=\mathbb {C} \) are valued in \( D \).
Hence by Theorem \ref{thrm:restricRtransformDvalued}, the \( D \)-valued cumulant
series of \( X \) are the same as the restriction of the \( B \)-valued cumulant
series; hence the only scalar-valued cumulant of \( X \) which is nonzero is
the second cumulant \( k_{11}^{X} \), so that the distribution of \( X \)
is the semicircle law.

Conversely, assume that the distribution of \( X \) is the semicircle law.
Let \( \eta (b)=k_{11}^{X}(b) \), \( b\in B \). Then we have
\begin{eqnarray*}
2\tau (\eta (1))^{2} & = & 2\tau (X^{2})^{2}\\
 & = & \tau (X^{4})=F\circ E(X^{4})\\
 & = & F(E(\eta (\eta (1))))+F(E(\eta (1)\eta (1)))\\
 & = & \tau (X\eta (1)X)+\tau (\eta (1)^{2})\\
 & = & \tau (\eta (1)XX)+\tau (\eta (1)^{2})\\
 & = & 2\tau (\eta (1)^{2}),
\end{eqnarray*}
so that
\[
\tau (\eta (1))=\tau (\eta (1)^{2})^{1/2}.\]
By the Cauchy-Schwartz inequality, we have that if \( \eta (1)\notin \mathbb {C} \),
\[
\tau (\eta (1))=\langle \eta (1),1\rangle <\Vert \eta (1)\Vert _{2}\cdot \Vert 1\Vert _{2}=\tau (\eta (1)^{2})^{1/2},\]
which is a contradiction. Hence \( \eta (1)\in \mathbb {C} \).
\end{proof}
We mention a corollary, which is of interest to random matrix theory. Let \( \sigma (x,y)=\sigma (y,x) \)
be a non-negative function on \( [0,1]^{2} \), having at most a finite number
of discontinuities in each vertical line. Let \( G(n) \) be an \( n\times n \)
random matrix with entries \( g_{ij} \), so that \( \{g_{ij}:i\leq j\} \)
are independent complex Gaussian random variables, \( g_{ij}=\overline{g_{ji}} \),
the expectation \( E(g_{ij})=0 \) and the variance \( E(|g_{ij}|^{2})=\frac{1}{n}\sigma (\frac{i}{n},\frac{j}{n}) \).
The matrices \( G(n) \) are called Gaussian Random Band Matrices. Let \( \mu _{n} \)
be the expected eigenvalue distribution of \( G(n) \), i.e.,
\[
\mu _{n}([a,b])=\frac{1}{n}\times \textrm{expected number of eigenvalues of }G(n)\, \textrm{in }[a,b].\]

\begin{cor}
The eigenvalue distribution measures \( \mu _{n} \) of the Gaussian Random
Band Matrices \( G(n) \) converge weakly to the semicircle law iff \( \int _{0}^{1}\sigma (x,y)dy \)
is a.e. a constant, independent of \( x \).
\end{cor}
The proof of this relies on a result from \cite{shlyakht:bandmatrix}, showing
that \( G(n) \) has limit eigenvalue distribution \( \mu  \), given as follows.
Let \( X \) be the \( L^{\infty }[0,1] \)-semicircular variable in an \( L^{\infty }[0,1] \)-probability
space \( (M,E:M\to L^{\infty }[0,1]) \), so that \( E(XfX)(x)=\int _{0}^{1}f(y)\sigma (x,y)dy \).
Let \( F:L^{\infty }[0,1]\to \mathbb {C} \) denote the linear functional \( F(f)=\int _{0}^{1}f(x)dx \),
and denote by \( \tau  \) the trace \( F\circ E \) on \( W^{*}(X,L^{\infty }[0,1]) \).
Then \( \mu  \) is the scalar-valued distribution of \( X \) with respect
to \( \tau  \), i.e.,
\[
\int t^{k}d\mu (t)=\tau (X^{k}).\]
It remains to apply Theorem \ref{thrm:characteriRestrictedSemicircular}, to
conclude that \( \mu  \) is a semicircle law iff \( E(X^{2})\in \mathbb {C} \),
i.e., \( \int _{0}^{1}\sigma (x,y)dy \) is a constant function of \( x \).

\subsection{A characterization of freeness.}

We are now ready to state the main theorem of this note. The following theorem
was earlier proved for \( B \)-valued semicircular variables in \cite{shlyakht:amalg},
and found many uses in operator algebra theory.

\begin{thm}   \label{theorem:main}
Let \( X_{1},\dots ,X_{n}\in M \).   
Assume that $F:B\to D$ satisfies
the faithfullness condition that if $b_1\in B$ and 
if $F(b_1 b_2)=0$ for all $b_2\in B$, then 
$b_1 = 0$.
Then \( X_{1},\dots ,X_{n} \) are free
from \( B \) with amalgamation over 
\( D \) iff their \( B \)-valued cumulant
series satisfies
\begin{equation}
\label{eqn:RtransformFactorizationCondition}
k^{X_{1},\dots ,X_{n}}_{B;i_{1},\dots ,i_{k}}(b_{1},\dots ,b_{k-1})=F(k^{X_{1},\dots ,X_{n}}_{B;i_{1},\dots ,i_{k}}(F(b_{1}),\dots ,F(b_{k-1})).
\end{equation}
for all \( b_{1},\ldots ,b_{k-1}\in B \). In short, \( k=F\circ k\circ F \).
(Here the cumulant series are computed in the \( B \)-probability space \( (M,E:M\to B) \)).
Equivalently,
\[
k^{X_{1},\dots ,X_{n}}_{B;i_{1},\dots ,i_{k}}(b_{1},\dots ,b_{k-1})=k^{X_{1},\dots ,X_{n}}_{D;i_{1},\dots ,i_{k}}(F(b_{1}),\dots ,F(b_{k-1})).\]

\end{thm}
We note that in the case that $M$ is a $C^*$-probability space,
the faithfulness assumption above is exactly the condition that the
GNS representation of $B$ with respect to the conditional expectation
$F$ is faithful.

Since \( X_{1},\dots ,X_{n} \) are free from \( B \) with amalgamation over
\( D \) iff the algebra \( N \) generated by \( X_{1},\dots ,X_{n} \) and
\( D \) is free from \( B \) over \( D \), the theorem above can be equivalently
stated as

\begin{thm}
Let \( N\subset M \) be a subalgebra of \( M \), containing \( D \). 
Assume that $F:B\to D$ satisfies
the faithfullness condition that if $b_1\in B$ and 
if $F(b_1 b_2)=0$ for all $b_2\in B$, then 
$b_1 = 0$.
Then
\( N \) is free from \( B \) over \( D \) iff
\[
\{n_{1}b_{1},n_{2}b_{2},\dots ,n_{k}\}_{E}=F\{n_{1}F(b_{1}),n_{2}F(b_{2}),\dots ,n_{k}\}_{E}\]
for all \( k \) and \( n_{1},\dots ,n_{k}\in N \), \( b_{1},\dots ,b_{k-1}\in B \). 
\end{thm}
\begin{proof}
We prove the theorem in the first formulation.

Assume that the condition (\ref{eqn:RtransformFactorizationCondition}) is satisfied
by the cumulant series of \( X_{1},\dots ,X_{n} \). Let \( Y_{1},\dots ,Y_{n} \)
be as in Section \ref{sec:canonical}. Since the freeness of \( X_{1},\dots ,X_{n} \)
from \( B \) with amalgamation over \( D \) is a condition on the \( B \)-moment
series of \( X_{1},\dots ,X_{n} \), and \( Y_{1},\dots ,Y_{n} \) have the
same \( B \)-moment series as \( X_{1},\dots ,X_{n} \), it is sufficient to
prove that \( Y_{1},\dots ,Y_{n}\in \mathcal{F} \) are free with amalgamation
over \( D \) from \( B \).

Since (\ref{eqn:RtransformFactorizationCondition}) is satisfied, \( \lambda _{j}^{0}=E_{B}(Y_{j})\in D \)
and hence \( Y_{1},\dots ,Y_{n} \) belong to the algebra \( \mathcal{L} \)
generated in \( \mathcal{F} \) by \( \lambda _{j}^{*} \) and \( \lambda _{p}^{q} \),
\( 1\leq j,p\leq n \), \( q\geq 1 \) and \( D \). Therefore, it is sufficient
to prove that \( \mathcal{L} \) is free from \( B \) with amalgamation over
\( D \).

Let \( w_{1},\dots ,w_{s}\in \mathcal{L} \), so that \( F\circ E_{B}(w_{j})=0 \),
and let \( b_{0},\dots ,b_{s}\in B \), so that \( F(b_{j})=0 \) (allowing
also \( b_{0} \) and/or \( b_{s} \) to be equal to \( 1 \)). We must prove
that
\[
F\circ E_{B}(b_{0}w_{1}b_{1}\cdots w_{s}b_{s})=0.\]
Note that the factorization condition (\ref{eqn:RtransformFactorizationCondition})
as well as the definition of the generators of \( \mathcal{L} \) (see \ref{relation:lambdas}
and \ref{relation:E}) imply that \( E_{B}|_{\mathcal{L}} \) has values in
\( D \). It follows that we may assume that \( E_{B}(w_{j})=0 \) (since \( F\circ E_{B}(w_{j})=E_{B}(w_{j})\in D \)).
By the definition of \( E_{B} \), its kernel is spanned by irreducible non-trivial
words in the generators \( \lambda _{j}^{*} \) and \( \lambda _{p}^{q} \).
Then
\[
W=b_{0}w_{1}\cdots w_{s}b_{s}\]
 is again a linear combination of words in the generators \( \lambda _{j}^{*} \)
and \( \lambda _{p}^{q} \). By linearity, we may reduce to the case that \( W \)
is a single word. If \( W \) is irreducible, it must be non-trivial (since
each \( w_{i} \) is non-trivial), hence \( E_{B}(W)=0 \), so that \( F\circ E_{B}(W)=0 \).
So assume that \( W \) is not irreducible. Since each \( w_{i} \) is irreducible,
this means that \( W \) contains a sub-word of the form
\begin{eqnarray*}
W & = & W_{1}\cdot d_{0}\lambda _{i_{1}}^{*}d_{1}\lambda ^{*}_{i_{2}}d_{2}\cdots \\
 &  & \quad d_{s}b_{1}d_{s+1}\lambda _{i_{s+2}}^{*}\cdots b_{j}\cdots \lambda _{i_{k}}^{*}d_{k}b_{r}d_{k+1}\lambda _{i_{k+1}}^{*}\\
 &  & \qquad \cdots d_{p}\lambda _{i_{p}}^{*}d_{p+1}\lambda _{j}^{p}\cdot W_{2}.
\end{eqnarray*}
Using the relation (\ref{relation:lambdas}) and the factorization condition
(\ref{eqn:RtransformFactorizationCondition}), we get that
\[
W=W_{1}d_{0}k_{i_{1},\dots ,i_{p},j}(d_{1},\dots ,F(d_{k}b_{r}d_{k+1}),\dots ,d_{p+1})W_{2}=0,\]
since \( F(d_{k}b_{r}d_{k+1})=d_{k}F(b_{r})d_{k+1}=0 \). Thus in any case,
\( F\circ E_{B}(W)=0 \).

We have therefore seen that the factorization condition implies freeness with
amalgamation.

To prove the other implication, assume that \( X_{1},\dots ,X_{n} \) are free
with amalgamation over \( D \) from \( B \). Let \( Z_{1},\dots ,Z_{n} \)
be \( B \)-valued random variables, so that
\[
k_{B;i_{1},\dots ,i_{k}}^{Z_{1},\dots ,Z_{n}}(b_{1},\dots ,b_{k-1})=F(k^{X_{1},\dots ,X_{n}}_{D,i_{1},\dots ,i_{k}}(F(b_{1}),\dots ,F(b_{k-1}))).\]
where \( k_{D} \) denote \( D \)-valued cumulants. (Note that the first occurrence
of \( F \) is actually redundant, as \( k^{X_{1},\dots ,X_{n}}_{D,i_{1},\dots ,i_{k}}(F(b_{1}),\dots ,F(b_{k-1})\in D \)).
Then by the first part of the proof, \( Z_{1},\dots ,Z_{n} \) are free from
\( B \) with amalgamation over \( D \). Moreover, by Theorem \ref{thrm:restricRtransformDvalued},
the \( D \)-valued distributions of \( Z_{1},\dots ,Z_{n} \) and \( X_{1},\dots ,X_{n} \)
are the same. By assumption, \( X_{1},\dots ,X_{n} \) are free with amalgamation
over \( D \) from \( B \). This freeness, together with the \( D \)-valued
distribution of \( X_{1},\dots ,X_{n} \), determines their \( B \)-valued
distribution.  Indeed, the freeness assumptions determine
$$ F\circ E  (b' b_0 X_{i_1} b_1 X_{i_2} \cdots b_{n-1} X_{i_n} b_n),\quad b',b_i\in B$$
which in view of the assumptions on $F$ determines $$
 E  (b_0 X_{i_1} b_1 X_{i_2} \cdots b_{n-1} X_{i_n} b_n),\quad b_i\in B.$$
 It follows that the \( B \)-valued distributions of \( X_{1},\dots ,X_{n} \)
and \( Z_{1},\dots ,Z_{n} \) coincide. Hence the \( B \)-valued cumulants
of \( X_{1},\dots ,X_{n} \) satisfy (\ref{eqn:RtransformFactorizationCondition}). 
\end{proof}
As an application, we have the following proposition (see e.g. \cite[Lemma 2.7]{shlyakht:cpentropy}):

\begin{prop}
Let \( N\subset M \) be a subalgebra. Let \( B\subset C\subset D\subset M \)
be subalgebras, and \( E_{C}:M\to C \), \( E_{B}:M\to B \), \( E_{D}:M\to D \)
be conditional expectations, so that \( E_{B}=E_{B}\circ E_{C} \), \( E_{C}=E_{C}\circ E_{D} \).  Assume that $E_C$ and $E_B$ and $E_D$ satisfy
the faithfullness assumptions of Theorem~\ref{theorem:main}.
Assume that \( N \) is free from \( C \) with amalgamation over \( B \),
and also free from \( D \) with amalgamation over \( C \). Then \( N \) is
free from \( D \) with amalgamation over \( B \).
\end{prop}
\begin{proof}
Since \( N \) is free from \( D \) with amalgamation over \( C \), we have
that for all \( n_{j}\in N \) and \( d_{j}\in D \),
\[
\{n_{1}d_{1},n_{2}d_{2},\dots ,n_{k}\}_{E_{D}}=E_{C}\{n_{1}E_{C}(d_{1}),\dots ,n_{k}\}_{E_{C}}.\]
Since \( N \) is free from \( C \) with amalgamation over \( D \), we get
similarly that for all \( c_{j}\in C \),
\[
\{n_{1}c_{1},n_{2}c_{2},\dots ,n_{k}\}_{E_{C}}=E_{B}\{n_{1}E_{B}(c_{1}),\dots ,n_{k}\}_{E_{B}}.\]
Applying this with \( c_{j}=E_{C}(d_{j}) \) and combining with the previous
equation gives
\[
\{n_{1}d_{1},n_{2}d_{2},\dots ,n_{k}\}_{E_{D}}=E_{B}\{n_{1}E_{B}(d_{1}),\dots ,n_{k}\}_{E_{B}},\]
since \( E_{B}=E_{B}\circ E_{C} \). Hence \( N \) is free from \( D \) with
amalgamation over \( B \).
\end{proof}
We mention as a corollary the following identity. Let \( B\subset C\subset D \)
be \( C^{*} \)-algebras, \( E^{N}_{B}:N\to B \), \( E^{C}_{B}:C\to B \) and
\( E^{D}_{C}:D\to C \) be conditional expectations, having faithful GNS representations.
Consider the reduced free product
\[
(((N,E^{N}_{B})*_{B}(C,E^{C}_{B})),E^{N}_{B}*\textrm{id})*_{C}(D,E^{D}_{C}),\]
where \( E^{N}_{B}*\textrm{id} \) denotes the canonical conditional expectation
from the free product \( (N,E^{N}_{B})*_{B}(C,E^{C}_{B}) \) onto \( C \).
Then
\begin{eqnarray*}
(((N,E^{N}_{B})*_{B}(C,E^{C}_{B})),E^{N}_{B}*\textrm{id})*_{C}(D,E^{D}_{C}) &  & \\
\cong (N,E^{N}_{B})*_{B}(D,E^{C}_{B}\circ E^{D}_{C}). &  & 
\end{eqnarray*}
To see this, it is sufficient to prove that \( N\subset (((N,E^{N}_{B})*_{B}(C,E^{C}_{B})),E^{N}_{B}*\textrm{id})*_{C}(D,E^{D}_{C}) \)
is free from \( D \) with amalgamation over \( B \), since both \( (((N,E^{N}_{B})*_{B}(C,E^{C}_{B})),E^{N}_{B}*\textrm{id})*_{C}(D,E^{D}_{C}) \)
and \( (N,E^{N}_{B})*_{B}(D,E^{C}_{B}\circ E^{D}_{C}) \) are generated by \( N \)
and \( D \) as \( C^{*} \)-algebras. But \( N \) is free from \( D \) over
\( C \), and from \( C \) over \( D \), by construction. Hence by the proposition
above, \( N \) is free from \( D \) over \( B \).

\providecommand{\bysame}{\leavevmode\hbox to3em{\hrulefill}\thinspace}


\begin{thebibliography}{1}

\bibitem{cumulants-nica}
A.~Nica, \emph{{$R$}-transforms of free joint distributions and non-crossing
  partitions}, Journal of Functional Analysis \textbf{135} (1996), 271--296.

\bibitem{shlyakht:bandmatrix}
D.~Shlyakhtenko, \emph{Random {Gaussian} band matrices and freeness with
  amalgamation}, Internat. Math. Res. Notices \textbf{20} (1996), 1013--1026.

\bibitem{shlyakht:amalg}
\bysame, \emph{Some applications of freeness with amalgamation}, J. reine
  angew. Math. \textbf{500} (1998), 191--212.

\bibitem{shlyakht:cpentropy}
\bysame, \emph{Free entropy with respect to a completely positive map}, Amer.
  J. Math. \textbf{122} (2000), no.~1, 45--81.

\bibitem{cumulants}
R.~Speicher, \emph{Multiplicative functions on the lattice of non-crossing
  partitions and free convolution}, Math. Annalen \textbf{298} (1994),
  193--206.

\bibitem{speicher:thesis}
R.~Speicher, \emph{Combinatorial theory of the free product with amalgamation
  and operator-valued free probability theory}, Mem. Amer. Math. Soc.
  \textbf{132} (1998), x+88.

\bibitem{dvv:amalg}
D.-V. Voiculescu, \emph{Operations on certain non-commutative operator-valued
  random variables}, Recent advances in operator algebras ({Orl\'eans, 1992}),
  no. 232, Ast{\'e}risque, 1995, pp.~243--275.

\bibitem{DVV:book}
D.-V. Voiculescu, K.~Dykema, and A.~Nica, \emph{Free random variables}, CRM
  monograph series, vol.~1, American Mathematical Society, 1992.

\end{thebibliography}
\end{document}